\documentclass[12pt]{article}

\usepackage{tabularx}
\usepackage{array}
\usepackage{booktabs}

\usepackage{amsmath, amssymb}
\usepackage[dvips]{graphicx}
\usepackage{amsthm, amscd}

\usepackage{bm}

\newcommand{\calA}{\mathcal{A}}
\newcommand{\calB}{\mathcal{B}}
\newcommand{\calC}{\mathcal{C}}
\newcommand{\calO}{\mathcal{O}}

\newcommand{\dS}{{\mathfrak S}}

\newcommand{\bbR}{{\mathbb R}}
\newcommand{\bbZ}{{\mathbb Z}}
\newcommand{\bbS}{{\mathbb S}}

\newcommand{\Ch}{{\rm \bf Ch}}
\newcommand{\rbS}{{\rm \bf S}}

\newtheorem{theorem}{Theorem}[section]
\newtheorem{prop}[theorem]{Proposition}

\newtheorem{lemma}[theorem]{Lemma}
\newtheorem{rem}[theorem]{Remark}

\def\Proof{\noindent{\sl Proof.}\qquad}

\begin{document}

\title{Arrangements stable under the Coxeter groups}
\author{
Hidehiko Kamiya 
\footnote
{
This work was partially supported by 
JSPS KAKENHI (22540134).
}
\\ 
{\it\footnotesize Graduate School of Economics, Nagoya University, 
Nagoya, 464-8601, Japan}
\\ 
Akimichi Takemura 
\footnote
{
This research was supported by JST CREST.
} 
\\ 
{\it\footnotesize Graduate School of Information Science and Technology,
University of Tokyo, 
Tokyo, 113-0033, Japan}
\\
Hiroaki Terao 
\footnote
{
This work was partially supported by 
JSPS KAKENHI (21340001).
}
\\ 
{\it\footnotesize Department of Mathematics, Hokkaido University, 
Sapporo, 060-0810, Japan}
} 
\date{October 2011}
\maketitle 

\begin{abstract}
Let $\calB$ be a real hyperplane arrangement 
which is stable under 
the action of 
a Coxeter group $W$. 
Then $W$ acts naturally 
on the set of chambers of $\calB$. 
We assume that $\calB$ is disjoint from the 
Coxeter arrangement $\calA=\calA(W)$ of $W$.
In this paper, we show that 
the 
$W$-orbits 
of the 
set of 
chambers of $\cal B$ are in one-to-one correspondence with 
the chambers of $\calC=\calA \cup \calB$ 
which are contained in an 
arbitrarily 
fixed chamber 
of $\cal A$.
From this fact, we find that the number of 
$W$-orbits of the set of chambers of $\calB$  is given by 
the number of chambers of $\calC$ divided by the order of $W$.
We will also study the set of 
chambers of $\calC$ which are 
contained in a chamber 
$b$ of $\calB$. 
We prove that the cardinality of this set is 
equal to the order of the isotropy subgroup 
$W_b$ of $b$. 
We illustrate these results with some examples, 
and 
solve 
an open problem in 
[H. Kamiya, A. Takemura, H. Terao, 
Ranking patterns of unfolding models of codimension one, 
Adv. in Appl. Math. 47 (2011) 379--400] 
by using our results.  
\end{abstract}





\section{Introduction} 
\label{sec:introduction}

Let $\calB$ be a real hyperplane arrangement 
which is stable under 
the action of 
a Coxeter group $W$. 
Then $W$ acts naturally 
on the set $\Ch(\calB)$ of chambers of $\calB$. 
We want to find the number of $W$-orbits of 
$\Ch(\calB)$. 
A particular case of this problem 
was considered in 
the authors' 
previous paper (Kamiya, Takemura and Terao \cite{ktt-aam}) and
the present 
paper is motivated by an open problem left in 
Section 6 of 
\cite{ktt-aam}. 
By the general results of the present paper, 
we give the affirmative answer to the open problem
in Theorem \ref{thm:affirmative}.
%

Suppose throughout that 
$\calB \cap \calA=\emptyset$, where 
$\calA=\calA(W)$ is the Coxeter arrangement 
of $W$. 
In this paper, we will show that the orbit space 
of $\Ch(\calB)$ is in one-to-one correspondence with 
the set of chambers $c$ of $\calC=\calA \cup \calB$ 
which are contained in 
$a$, $\{ c \in \Ch(\calC) \mid c \subseteq a \}$, 
where $a \in \Ch(\calA)$ is an arbitrary chamber 
of $\calA$.  
From this fact, we find that the number of 
$W$-orbits of $\Ch(\calB)$ is given by 
$|\Ch(\calC)|/|W|$. 

On the other hand, we will also study the set of 
chambers $c \in \Ch(\calC)$ which are 
contained in a chamber 
$b \in \Ch(\calB)$ of $\calB$, 
$\{ c \in \Ch(\calC) \mid c \subseteq b\}$.  
We will prove that the cardinality of this set is 
equal to the order of the isotropy subgroup 
$W_b$ 
of $b$. 
Moreover, we will investigate the structure of $W_b$. 


Kamiya, Takemura and Terao \cite{ktt-aam} 
tried to find 
the number of ``inequivalent ranking patterns 
generated by unfolding models of codimension one'' 
in psychometrics, and obtained an upper bound for this 
number. 
It was left open to determine whether this upper bound 
is actually the exact number. 
The problem boils down to 
proving (or disproving) that 
the orbit space of the chambers of 
the restricted all-subset arrangement 
(\cite{ktt-aam})  
$\calB$ 
under the action of the symmetric group 
$\dS_m$ is in one-to-one correspondence with 
$\{ c \in \Ch(\calA(\dS_m)\cup \calB) \mid 
c \subseteq a\}$ 
for a chamber $a \in \Ch(\calA(\dS_m))$ of the braid 
arrangement $\calA(\dS_m)$. 
The results of the present paper establish the
one-to-one correspondence.

The paper is organized as follows. 
In Section \ref{sec:main}, we verify our main results. 
Next, in Section \ref{sec:example}, 
we illustrate our general results with 
five examples, some of which 
are 
taken from the authors' 
previous studies of unfolding models in psychometrics 
(\cite{kott}, \cite{ktt-aam}). 
In Section \ref{sec:example}, 
we also solve the open problem of \cite{ktt-aam} 
(Theorem \ref{thm:affirmative}) 
using our general results in Section \ref{sec:main} 
applied to one of our examples.   

\section{Main results}
\label{sec:main}

In this section, we 
state and prove 
our main results. 

Let $V$ be 
a 
Euclidean space.
Consider a Coxeter group $W$ acting on $V$. 
Then the Coxeter arrangement 
$\calA=\calA(W)$ is the set of all reflecting hyperplanes of $W$. 
Suppose that $\calB$ is a hyperplane arrangement which is 
stable under the natural action of $W$. 
We assume 
$\calA \cap \calB=\emptyset$ and define
\[
\calC:=\calA \cup \calB. 
\]

Let $\Ch(\calA)$, $\Ch(\calB)$ and
$\Ch(\calC)$ 
denote the set of chambers of $\calA$, $\calB$ 
and $\calC$, respectively.
Define 
\begin{align*}
\varphi_{\calA}:\Ch(\calC)\to \Ch(\calA),
\quad 
\varphi_{\calB}:\Ch({\calC})\to \Ch({\calB}) 
\end{align*}
by
\begin{align*}
\varphi_{\calA} (c)
:=\text{the chamber of $\calA$ containing $c$},\\
\varphi_{\calB}(c)
:=\text{the chamber of $\calB$ containing $c$}
\end{align*}
for $c\in \Ch(\calC)$. 
Note that the Coxeter group $W$ naturally acts on $\Ch(\calA)$,
$\Ch(\calB)$ and $\Ch(\calC)$.
\begin{lemma}\label{lemma:1}
$\varphi_{\calA}$ and $\varphi_{\calB}$ are
both $W$-equivariant,
i.e.,
\begin{align*}
\varphi_{\calA}(wc)=w(\varphi_{\calA}(c)),
\quad 
\varphi_{\calB}(wc)=w(\varphi_{\calB}(c))
\end{align*}
for any $w \in W$ and $c \in \Ch(\calC)$.
\end{lemma}

The proof is easy and omitted. 

The following result is classical
(see, e.g., \cite[Ch. V, \S 3. 2. Theorem 1 
(iii)
]{bou}):
\begin{theorem}\label{theorem:2}
The group $W$ acts on $\Ch(\calA)$ effectively and transitively.
In particular, $|W|=|\Ch(\calA)|$.
\end{theorem}
Using Theorem \ref{theorem:2}, we can prove the following 
lemma. 
\begin{lemma}\label{lemma:3}
The group $W$ acts on $\Ch(\calC)$ effectively. In particular,
each $W$-orbit of $\Ch(\calC)$ is of size $|W|$.
\end{lemma}

\Proof 
If $wc=c$ for $w\in W$ and $c\in\Ch(\calC)$,
then we have 
$
\varphi_{\calA}(c)
=
w
\varphi_{\calA}(c),
$   
which implies $w=1$ by Theorem \ref{theorem:2}.
\qed \\ 

For $b\in \Ch(\calB)$, define 
the isotropy subgroup
$W_b:=\{w\in W \mid wb=b\}$.
Then we have 
the next lemma. 
\begin{lemma}\label{lemma:4}
For $b \in \Ch(\calB)$, 
the group $W_b$ acts on 
$\varphi_{\calB}^{-1}(b)$ effectively and transitively.
\end{lemma}

\Proof 
The effective
part follows from 
Lemma \ref{lemma:3}, 
so let us prove the transitivity. 
Let $c_1, c_2 \in \varphi_{\calB}^{-1}(b)$.
Define
$$\calA(c_1, c_2):=\{H \in \calA 
\mid 
c_1
{\text{~and~}} c_2 {\text{ are on different sides of~}} H\}.$$
Let us prove that there exists $w\in W$ such that
$wc_1=c_2$ by an induction on $|\calA(c_1, c_2)|$.
When $|\calA(c_1, c_2)|=0$, 
we have 
$\calA(c_1, c_2)=\emptyset$ and
$c_1=c_2$. Thus we may choose $w=1$. 
If
$\calA(c_1, c_2)$ is non-empty, then there exists
$H_1\in \calA(c_1, c_2)$ such that $H_1$ 
contains a wall of 
$c_{1}$. Let $s_1$ denote the reflection with respect 
to $H_1$. Then
$$\calA(s_1 c_1, c_2)=\calA(c_1, c_2)
\setminus\{H_1\}.$$
By the induction assumption, there exists $w_{1} \in W$ 
with $w_{1} s_{1} c_{1} = c_{2} $.
Set $w:=w_{1} s_{1} $.
Then
$w c_{1} = c_{2} $
 and $c_{2} = (w c_{1}) \cap c_{2} \subseteq (wb)\cap b$, 
which implies that $(wb)\cap b$ is not empty.
Thus 
$wb=b$ and $w\in W_{b} $. 
\qed \\ 

The following lemma states that 
the $W$-orbits of $\Ch(\calC)$ and those of 
$\Ch(\calB)$ are in one-to-one correspondence. 
\begin{lemma}\label{lemma:5}
The map $\varphi_{\calB} : \Ch(\calC) \rightarrow \Ch(\calB)$
induces a bijection
from 
the set of $W$-orbits of $\Ch(\calC)$
to
the set of $W$-orbits of $\Ch(\calB)$.
\end{lemma}

\Proof 
For $b\in \Ch(\calB)$ and $c\in \Ch(\calC)$, we denote
the $W$-orbit of $b$ and
the $W$-orbit of $c$ 
by $\calO(b)$
and  
by $\calO(c)$, 
respectively.  
It is easy to see that 
\[
\varphi_{\calB}(\calO(c))
=\calO(\varphi_{\calB}(c)), 
\quad 
c \in \Ch(\calC), 
\]
by Lemma \ref{lemma:1}.
Thus $\varphi_{\calB}$ induces a map
from 
the set of $W$-orbits of $\Ch(\calC)$
to
the set of $W$-orbits of $\Ch(\calB)$.
We will show the map is bijective.

{\em Surjectivity}: 
Let $\calO(b)$ be an arbitrary 
orbit of $\Ch(\calB)$ 
with a representative point $b \in \Ch(\calB)$. 
Take an arbitrary $c\in\varphi_{\calB}^{-1}(b)$.    
Then   
$$
\varphi_{\calB}  (\calO(c))
=
\calO(\varphi_{\calB}(c))
=
\calO(b),
$$ 
which shows the surjectivity.

{\em Injectivity}: 
Suppose 
$
\varphi_{\calB}  (\calO(c_{1} ))
=
\varphi_{\calB}  (\calO(c_{2} )) 
\,\, 
(c_{1}, c_{2} \in \Ch(\calC))$.
Set $b_{i} :=
\varphi_{\calB}  (c_{i})$ 
for $i=1,2$.  
We have 
$$
\calO(b_{1} )
=
\calO(\varphi_{\calB}  (c_{1}))
=
\varphi_{\calB}(\calO(c_{1}))
=
\varphi_{\calB}(\calO(c_{2}))
=
\calO(\varphi_{\calB}  (c_{2}))
=
\calO(b_{2} ), 
$$
so we can pick $w\in W$ such that 
$w b_{2} = b_{1} $.
  Then
\[
\varphi_{\calB}(w c_{2})
=
w
(\varphi_{\calB}(c_{2}))
=
w
b_{2}
=
b_{1}.
\]
Therefore, 
both $c_{1} $ and
$w c_{2} $ lie in
$\varphi_{\calB}^{-1} (b_{1} ) $.
By Lemma \ref{lemma:4}, we have
$
\calO(c_{1} )
=
\calO(w c_{2} )
=
\calO(c_{2} ).
$
\qed \\ 

We are now in a position to state the main results 
of this paper. 

\begin{theorem}\label{theorem:6}

The cardinalities of 
$\varphi_{\calA}^{-1} (a), \ 
\varphi_{\calB}^{-1} (b)$ 
for $a \in \Ch(\calA), \ b \in\Ch(\calB)$ 
are given as follows:    

\begin{enumerate}
\item 
For $a\in\Ch(\calA)$, 
we have 
\begin{align*} 
|\varphi_{\calA}^{-1} (a) |
=
\frac{|\Ch(\calC)|}{|\Ch(\calA)|}  
=
\frac{|\Ch(\calC)|}{|W|}
&=
|
\{
W\text{-orbits of~} \Ch(\calC)
\}  
|\\
&=
|
\{
W\text{-orbits of~} \Ch(\calB)
\}  
|.
\end{align*} 
\item 
For $b\in\Ch(\calB)$, 
we have 
$
|\varphi_{\calB}^{-1} (b) |
=
|W_{b}|.
$ 
\end{enumerate}
\end{theorem}

\Proof 
Part 2 follows 
from Lemma \ref{lemma:4}, so 
we will prove Part 1. 
Since the map $\varphi_{\calA} : \Ch(\calC) \rightarrow \Ch(\calA)$
is $W$-equivariant (Lemma \ref{lemma:1}), 
we have 
for each $w \in W$ 
a bijection
\[
\varphi_{\calA}^{-1}(a)
\rightarrow  
\varphi_{\calA}^{-1}(wa)
\]
sending $c\in \varphi_{\calA}^{-1}(a)
$ to
$wc$. 
Thus every fiber of $\varphi_{\calA}
$  has the same cardinality because
$W$ acts transitively on $\Ch(\calA)$
(Theorem \ref{theorem:2}).
The cardinality is equal to
$$
\frac{|\Ch(\calC)|}{|\Ch(\calA)|}  
=
\frac{|\Ch(\calC)|}{|W|}.
$$  
By Lemma \ref{lemma:3}, we have 
\[
|
\{
W\text{-orbits of~} \Ch(\calC)
\}  
|
=
\frac{|\Ch(\calC)|}{|W|}. 
\]
Finally, Lemma \ref{lemma:5} proves the last equality.
\qed \\

By Part 2 of Theorem \ref{theorem:6}, 
we can write $|\Ch(\calC)|$ as 
\begin{eqnarray} 
\label{eq:Fix}
|\Ch(\calC)|
&=& \sum_{b\in \Ch(\calB)}|\varphi_{\calB}^{-1}(b)|
=\sum_{b\in \Ch(\calB)}|W_b| \nonumber \\ 
&=& \sum_{w\in W}{\rm Fix}(w,\Ch(\calB)), 
\end{eqnarray} 
where ${\rm Fix}(w,\Ch(\calB))$ denotes the number of 
elements of $\Ch(\calB)$ fixed by $w$. 
In the case of the Catalan arrangement, 
\eqref{eq:Fix} is stated in 
\cite[p. 561]{ps-00}. 

Next, let $x \in V\setminus \bigcup_{H \in \calB}H$.
Let $b \in \Ch(\calB)$ denote the unique chamber 
that contains $x$. 
Define the average $z(x)$ 
of $x$ over the action of $W_b$: 
\[
z(x)=
\frac{1}{|W_b|}
\sum_{w \in W_b}wx. 
\]
Then it is easily seen that $z(x)$ lies in $b$ because 
of the convexity of $b$, and that the map $z$ is 
$W$-equivariant.  
Concerning the structure of $W_b$, 
we have the next proposition. 

\begin{prop}\label{prop:W_b=W_z}
The following statements hold true: 
\begin{enumerate}
\item 
For any $b \in \Ch(\calB)$, 
the set
$
\{ z(x) \mid x \in b\} 
$
is equal to
the set of all $W_b$-invariant points of $b$. 
\item 
For any $x \in V\setminus \bigcup_{H\in \calB}H$, 
the isotropy subgroup
$W_{z(x)}$ 
of $z(x)$ 
is equal to $W_{b}$, 
where $b \in \Ch(\calB)$ is the unique chamber 
that contains $x$. 
In particular, 
$W_{z(x)}$ depends only on the chamber
$b\in \Ch(\calB)$ containing
$x$.

\end{enumerate}
\end{prop}

\Proof 


1.  
By the linearity of the action of $W$ on $V$,
the average $z(x)\in b$ of $x\in b$  
is $W_{b}$-invariant: 
$wz(x)=z(x), \ w \in W_b$. 
Conversely, the average of any $W_{b}$-invariant point 
of $b$ 
is the point itself. 
  
2.  
Assume
$w \in W_{z(x)}$.
Then $z(x)=wz(x)  \in b \cap wb$, which implies 
$b \cap wb \ne
\emptyset$. 
Since $b$ and $wb$ are both chambers, they coincide:
$wb=b$. 
Thus 
$w\in W_{b} $
and 
we obtain 
$W_{z(x)  } \subseteq W_{b}. $ 
We also have the reverse inclusion because 
the average 
$z(x)$ is $W_{b}$-invariant by 
the statement 1. 
\qed \\ 

%

When $W$ is the symmetric group 
$\dS_m=\dS_{\{ 1,\ldots,m\}}$, 
we have the following obvious fact 
(Proposition \ref{prop:W_b}). 
Define 
\[
H_0:=
\{ x=(x_1,\ldots,x_m)^T \in \bbR^m 
\mid x_1+\cdots +x_m=0\}. 
\]
The group $W=\dS_m$ acts on $V=\bbR^m$ or $V=H_0$ 
by permuting coordinates.  
When $W=\dS_m$ and $V=\bbR^m$ or $V=H_0$, 
we 
agree 
that this action is considered. 

\begin{prop}\label{prop:W_b}
Let $W=\dS_m$ and 
$V=\bbR^m$ 
or $V=H_0$. 
Then we have 
\[
W_b = 
\dS_{k_1} \times \dS_{k_2} 
\times \cdots \times \dS_{k_{\ell}}, 
\quad b \in \Ch(\calB),  
\]
where $k_1,\ldots,k_{\ell}$ 
($k_1+\cdots +k_{\ell}=m, \ 1\le \ell \le m$) 
are defined by 
\[
z_{\sigma(1)}= \cdots =z_{\sigma(k_1)}
>z_{\sigma(k_1+1)}= \cdots =z_{\sigma(k_1+k_2)}
>\cdots 
>z_{\sigma(k_1+\cdots +k_{\ell -1}+1)}
= \cdots =z_{\sigma(m)}
\]
for 
$z=(z_1,\ldots,z_m)^T=z(x), \ x \in b$,  
and 
a permutation 
$\sigma \in \dS_{\{ 1,\ldots,m\}}$. 
\end{prop}

\begin{rem}
\label{rem:conjugacy-injection}
Consider the map $b \mapsto W_b$ 
from $\Ch(\calB)$ to the set of 
subgroups of $W$.
Since $W_{wb}=w W_b w^{-1}, \ w\in W$, 
this induces a map $\tau$ 
from the set of $W$-orbits of $\Ch(\calB)$
to the set of conjugacy classes of subgroups of $W$: 
\begin{equation}
\label{eq:def-tau}
\tau(\calO(b))=[W_b], 
\quad b \in \Ch(\calB), 
\end{equation}
where $[W_b]:=\{ w W_b w^{-1} \mid w \in W\}$. 
This map $\tau$ is not injective in general.  
See Remarks \ref{rem:non-unique-example-semiorder} 
and \ref{rem:non-unique-example}. 
\end{rem}

\section{Examples}
\label{sec:example}

In this section, we examine 
five examples. 
The first example (Subsection \ref{subsec:Catalan}) is 
the Catalan arrangement, which has been 
well studied 
as a deformation of the braid arrangement. 
The Catalan arrangement is related to the semiorder 
introduced by Luce \cite{luce-56} in economics and 
mathematical psychology as a preference order that 
accounts for 
intransitive indifference. 
The next three examples are taken from problems 
in psychometrics---the arrangements in Subsections 
\ref{subsec:A+r-allsubset} and 
\ref{subsec:A+u-allsubset} 
(the braid arrangement in conjunction with 
the all-subset arrangement) 
appear 
naturally 
in the study of ranking patterns of 
unfolding models of codimension one 
(Kamiya, Takemura and Terao \cite{ktt-aam}), 
while the mid-hyperplane arrangement in 
Subsection \ref{subsec:mid-hyper} is needed 
in 
examining 
ranking patterns of 
unidimensional unfolding models 
(Kamiya, Orlik, Takemura and Terao \cite{kott}). 
In all four examples in Subsections \ref{subsec:Catalan}--\ref{subsec:mid-hyper}, 
the Coxeter group $W$ is 
of type $A_{m-1}$. 
In Subsection \ref{subsec:signed-allsubset}, 
we provide 
an illustration 
with the Coxeter group of type $B_m$. 
We also solve the open problem of \cite{ktt-aam} 
in Subsection \ref{subsec:A+r-allsubset}. 

\subsection{Catalan arrangement} 
\label{subsec:Catalan} 




Let $W=\dS_m$ and $V=H_0$. 
Then $\calA=\calA(W)$ is the braid arrangement 
in $H_0$, consisting of the hyperplanes defined by 
$x_i=x_j, \ 1\le i<j\le m$. 
All the $|\calA|
=m(m-1)/2$
hyperplanes form one orbit under the action of 
$W$ on $\calA$.  

Let 
\[
\calB=
\{ H_{ij} \mid 
i\ne j, \ 1\le i \le m, \ 1\le j\le m 
\}, 
\]
where 
\[
H_{ij}:=\{ x=(x_1,\ldots,x_m)^T\in H_0 
\mid 
x_i=x_j+1 
\}. 
\]
That is, $\calB$ is an essentialization of 
the semiorder arrangement. 
We have $|\calB|=m(m-1)$, 
and the action of $W$ on $\calB$ is transitive. 
 
With these $\calA$ and $\calB$, the union 
$\calC=\calA \cup \calB$ is 
an essentialization of 
the Catalan arrangement. 

Now, as a chamber of $\calA$, let us take 
$a \in \Ch (\calA)$ defined by $x_1>\cdots >x_m$: 
\begin{equation}
\label{eq:a:x1>xm} 
a: x_1>\cdots >x_m. 
\end{equation} 
For this $a$, 
we have always $x_i-x_j<1 \ (i>j)$ for 
$(x_1,\ldots,x_m)^T\in a$, 
so the elements of 
$\varphi_{\calA}^{-1}(a)$ are determined by 
the sets of pairs $i,j \ (i<j)$ such that 
$x_i-x_j<1$, or equivalently 
the sets of maximal intervals 
$[i,j]:=\{ i, i+1, \ldots,j\} 
\ (i<j)$ 
such that $x_i-x_j<1$, 
for $(x_1,\ldots,x_m)^T\in a$.   
It is well known 
(\cite[(1.1)]{wf-57}, 
\cite{dk-68}) 
that the number $|\varphi_{\calA}^{-1}(a)|$ is 
equal to the Catalan number 
\[
C_m
:=\frac{1}{m+1}\binom{2m}{m}. 
\]
Thus  
\[
|\Ch (\calC)|=|\varphi_{\calA}^{-1}(a)|\times |W|
=m! \, C_m
=
(2m)_{m-1} := 
(2m)(2m-1)\cdots (m+2). 
\]  
Besides, 
from $|\Ch (\calC)|=m! \, C_m$,  
the characteristic polynomial $\chi(\calC, t)$ 
of $\calC$ (\cite[Definition 2.52]{ort}) 
can be calculated as  
\[
\chi(\calC, t)=(t-m-1)(t-m-2)\cdots (t-2m+1) 
\]
(\cite[Theorem 5.18]{sta-07}). 


The $W$-orbits of $\Ch(\calC)$ are given as follows. 
Let 
$c_i \in \Ch(\calC), \ i=1,\ldots,C_m$, 
be the chambers of $\calC$ 
which are 
contained in 
$a$ in \eqref{eq:a:x1>xm}: 
$\varphi_{\calA}^{-1}(a)= 
\{ c_i \mid i=1,\ldots,C_m\}$. 
As was mentioned earlier, each $c_i, \ i=1,\ldots,C_m$, can be 
indexed by the set of maximal intervals  
$[i,j]\subseteq \{ 1,\ldots,m\} \ (i<j)$ such that 
$x_i-x_j<1$.
The set 
$\varphi_{\calA}^{-1}(a)
=\{ c_i \mid i=1,\ldots,C_m\}
\subset \Ch(\calC)$ is 
a complete set of representatives of the
$W$-orbits of $\Ch(\calC)$, 
i.e., $\Ch(\calC)$ has exactly $C_m$ orbits 
$\calO(c_i), \ i=1,\ldots,C_m$,  
under the action of $W$. 

Next, consider the chambers of $\calB$. 
The elements of $\Ch(\calB)$ are 
in one-to-one correspondence with the set of semiorders 
on $\{ 1,\ldots,m\}$. 
Recall that a partial order 
$\succeq$ on $\{ 1,\ldots,m\}$ is 
called a {\it semiorder} (or a {\it unit interval order}) 
if and only if 
the poset $(\{ 1,\ldots,m\}, \succeq)$ 
contains no induced subposet 
isomorphic to 
$\bm{2}+\bm{2}$ or $\bm{3}+\bm{1}$, 
where $\bm{i}+\bm{j}$ stands for 
the disjoint union of an $i$-element chain 
and a $j$-element chain. 
The Scott-Suppes Theorem (\cite{ss-58}) states that 
a partial order $\succeq$ on $\{ 1,\ldots,m\}$ is a semiorder 
if and only if there exist $x_1,\ldots,x_m \in \bbR$ such that 
$i \succeq j \iff i \succ j \text{ or } i=j 
\ (i,j\in \{ 1,\ldots,m\})$, 
where 
\begin{equation}
\label{eq:i-succ-j} 
i \succ j \iff x_i>x_j+1. 
\end{equation} 
Now, a bijection from $\Ch (\calB)$ to the set $\rbS(m)$ 
of semiorders on $\{ 1,\ldots,m\}$ is given by 
$\Ch (\calB)\ni b \mapsto \ \succeq \ \in \rbS(m)$, where 
$\succeq$ is the partial order determined by \eqref{eq:i-succ-j} with 
an arbitrary 
$(x_1,\ldots,x_m)^T \in b$. 
(The surjectivity of this map follows from the 
``Distinguishing Property'' of $\rbS(m)$ 
(\cite[Lemma 5.1]{ovc-05}, \cite[Lemma 7.35]{ovc-11}).) 

The number of chambers of $\calB$ is obtained as follows. 
Let us denote (only in this paragraph) 
our $\calB$ and $\calC$ as $\calB_m$ and $\calC_m$, respectively. 
Postnikov and Stanley \cite[Lemma 7.6]{ps-00} 
proved that 
${\rm Fix}(w,\Ch(\calB_m))=|\Ch(\calB_k)|$, where $k$ is the number of 
cycles of $w\in W=\dS_m$. 
Thus \eqref{eq:Fix} is written as 
\begin{equation}
\label{eq:|Ch(C)|=sum-c|Ch(B)|} 
|\Ch(\calC_m)|=\sum_{k=1}^{m}c(m,k)|\Ch(\calB_k)|, 
\quad m\ge 1, 
\end{equation} 
in this case 
(\cite[Theorem 7.1]{ps-00}, \cite[Solution to Exercise 6.30]{sta-ec2}), 
where $c(m,k)$ is the signless Stirling number of the first kind. 
Because of 
\cite[Proposition 1.4.1]{sta-ec1}, 
equation \eqref{eq:|Ch(C)|=sum-c|Ch(B)|} is equivalent to 
\begin{equation}  
\label{eq:|Ch(B)|=sum(-1)S|Ch(C)|=} 
|\Ch(\calB_m)| 
=\sum_{k=1}^{m}(-1)^{m-k}S(m,k)|\Ch(\calC_k)| 
=\sum_{k=1}^{m}(-1)^{m-k}S(m,k)(2k)_{k-1}, 
\quad m\ge 1, 
\end{equation} 
where $S(m,k)$ is the Stirling number of the second kind. 
The number of semiorders, $|\rbS(m)|=|\Ch(\calB_m)|$, 
in \eqref{eq:|Ch(B)|=sum(-1)S|Ch(C)|=} 
was first obtained by 
Chandon, Lemaire and Pouget 
\cite[Proposition 12]{clp-78}. 
In addition, 
by \eqref{eq:|Ch(B)|=sum(-1)S|Ch(C)|=} 
and $\sum_{m=1}^{\infty}S(m,k)x^m/(m!)=(e^x-1)^k/(k!), \ k\ge 1$ 
(\cite[p. 174]{maz-10}), 
we obtain 
\begin{eqnarray} 
\label{eq:DGF-|Ch(B)|} 
\sum_{m=1}^{\infty}|\Ch(\calB_m)| \, \frac{x^m}{m!}
&=& \sum_{k=1}^{\infty}|\Ch(\calC_k)| \, \frac{(1-e^{-x})^k}{k!} \nonumber \\ 
&=& C(1-e^{-x}) \nonumber \\ 
&=& 1\cdot x+3\cdot \frac{x^2}{2!}
+19\cdot \frac{x^3}{3!}
+183\cdot \frac{x^4}{4!} \\ 
&& \qquad +2371\cdot \frac{x^5}{5!}
+38703\cdot \frac{x^6}{6!} +763099\cdot \frac{x^7}{7!}
+\cdots \nonumber 
\end{eqnarray} 
(\cite[Theorem 7.1]{ps-00}, 
\cite[Corollary 5.12]{sta-07}, 
\cite[Exercise 6.30]{sta-ec2}; 
cf. \cite[Table on p. 79]{clp-78}), 
where 
\[
C(t):=
\sum_{k=1}^{\infty}C_k t^k
=\frac{1-\sqrt{1-4t}}{2t}-1
=t+2t^2+5t^3
+\cdots 
\]  
(\cite[p. 178]{sta-ec2}). 
By \eqref{eq:|Ch(B)|=sum(-1)S|Ch(C)|=} and 
$\sum_{m=1}^{\infty}S(m,k)x^m
=\prod_{j=1}^k \{ x/(1-jx)\}, \ k \ge 1$ 
(\cite[Theorem 4.3.1]{maz-10}), 
we can also get 
\[
\sum_{m=1}^{\infty}|\Ch(\calB_m)| \, x^m 
=\sum_{k=1}^{\infty} \frac{(2k)_{k-1} \, x^k}{(1+x)(1+2x)\cdots (1+kx)}. 
\]

As for the $W$-orbits of $\Ch(\calB)$, 
we have 
$\{ W\text{-orbits of~} \Ch(\calB)\}
=\{ \calO(b_1), \ldots, \calO(b_{C_m})\}$ 
with cardinality $C_m$, 
where $b_i:=\varphi_{\calB}(c_i)\in \Ch(\calB), \ i=1,\ldots,C_m$. 

\ 

Now, let us investigate the case $m=3$. 

\begin{figure}[htbp]
 \begin{center}
\includegraphics*[width=.6\textwidth]{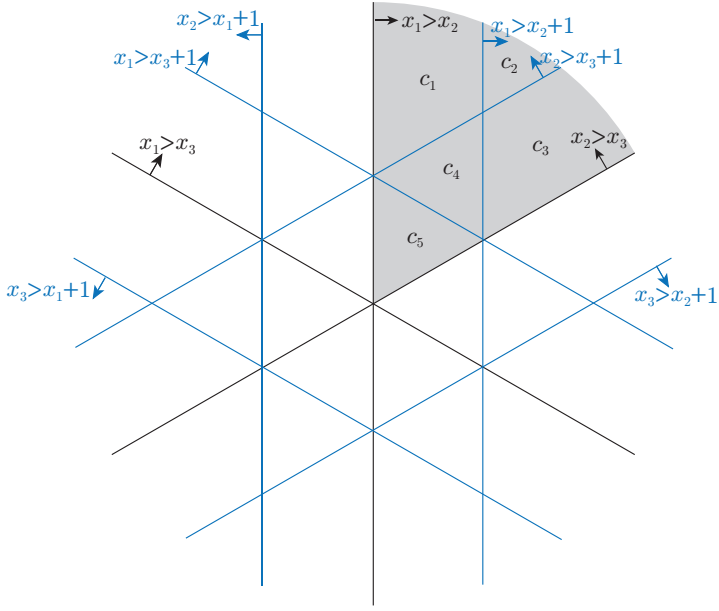}
\caption{Essentialization of Catalan arrangement.}
\label{fig:catalan-m=3}
 \end{center}
\end{figure}

The arrangement $\calA$ consists of three lines 
in $V=H_0, \ \dim H_0=2$, 
each of which is 
defined by one of the following equations:   
\begin{equation}
\label{eq:braid-with-m=3}
x_1=x_2, \quad x_1=x_3, \quad x_2=x_3, 
\end{equation}
and $\calB$ comprises six lines:  
\begin{equation*}
\label{eq:semiorder-m=3}
x_1=x_2\pm 1, \quad x_1=x_3\pm 1, \quad x_2=x_3\pm 1. 
\end{equation*}
Figure \ref{fig:catalan-m=3} 
displays $\calA$ and $\calB$ in $V=H_0$. 

Let $a \in \Ch (\calA)$ be as in \eqref{eq:a:x1>xm}: 
$x_1>x_2>x_3$ 
($a$ is shaded in 
Figure \ref{fig:catalan-m=3}). 
This chamber $a$ of $\calA$ contains exactly 
five $(=C_3)$ 
chambers $c_1,\ldots,c_5$ of $\calC$: 
$\varphi_{\calA}^{-1}(a)=\{ c_1,\ldots,c_5\}$. 
These five chambers $c_1,\ldots,c_5$ are indexed by 
the following sets of maximal intervals $[i,j]$ such 
that $x_i-x_j<1$: 
\begin{eqnarray*}
&& c_1: \{ [1,2]\}, \\ 
&& c_2: \, \emptyset, \\ 
&& c_3: \{ [2,3]\}, \\  
&& c_4: \{ [1,2], \, [2,3]\}, \\   
&& c_5: \{ [1,3]\}. 
\end{eqnarray*} 
Note that $c_1$ and $c_3$ are obtained 
from each other by changing $(x_1,x_2,x_3)$ to 
$(-x_3,-x_2,-x_1)$. 
Since $a \in \Ch(\calA)$ consists of 
$|\varphi_{\calA}^{-1}(a)|=5$ chambers 
$c_1,\ldots,c_5$ of $\calC$, 
we have $|\Ch(\calC)|
=|\varphi_{\calA}^{-1}(a)| \times |W|=5\times 3!=30$, 
and $\Ch(\calC)$ has exactly five $W$-orbits 
$\calO(c), \ldots,\calO(c_5)$. 

\begin{figure}[htbp]
 \begin{center}
\includegraphics*[width=.7\textwidth]{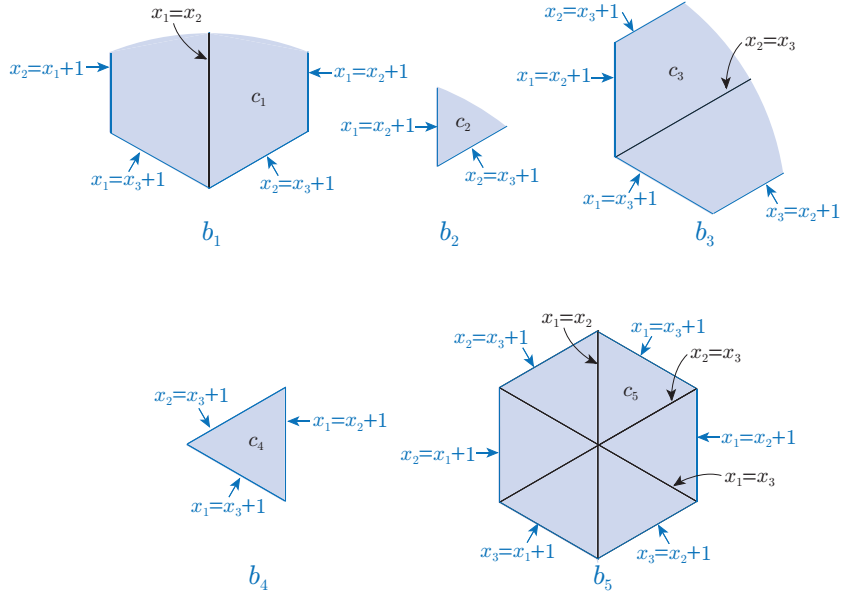}
\caption{Essentialization of Catalan arrangement. 
$b_1,\ldots, b_5$.}
\label{fig:catalan-m=3-b}
 \end{center}
\end{figure}

In $c_1: \{ [1,2]\}$, we have 
$x_1-x_2<1, \ x_1-x_3>1$ and $x_2-x_3>1$, 
so the semiorder corresponding to $b_1=\varphi_{\calB}(c_1)$ 
is $1 \succ 3, \ 2 \succ 3$. 
(Recall that for the $a$ in \eqref{eq:a:x1>xm}, and hence for 
any $c \in \Ch(\calC)$ contained in $a$, 
we have never $i \succ j$ for $i>j$.) 
By similar arguments, we can also obtain the semiorders corresponding to 
$b_i=\varphi_{\calB}(c_i)$ for $i=2,\ldots,5$:  
\begin{eqnarray}
\label{eq:b1b5-semiorder}
&& b_1: 1 \succ 3, \ 2 \succ 3, \nonumber \\ 
&& b_2: 1 \succ 2 \succ 3, \nonumber \\  
&& b_3: 1 \succ 2, \ 1 \succ 3, \\ 
&& b_4: 1 \succ 3, \nonumber \\ 
&& b_5: \text{none}. \nonumber 
\end{eqnarray} 
See Figure \ref{fig:catalan-m=3-b}. 
Again, $b_1$ and $b_3$ are obtained 
from each other by the above-mentioned rule. 
The chamber $b_1 \in \Ch(\calB)$ is divided 
by the line of $\calA$ defined by $x_1=x_2$ 
into two chambers of $\calC$, 
$|\varphi_{\calB}^{-1}(b_1)|=2$. 
The $W_{b_1}$-invariant points $z$ of $b_1$ are 
$z=d(1,1,-2)^T, \ d>1/3$, 
so we have 
$W_{b_1}=\dS_{\{ 1,2\}}$. 
In a similar manner, 
$b_2$ is not divided by any line in $\calA$;  
$b_3$ is divided by the line 
$x_2=x_3$ into two;  
$b_4$ is not divided by any line; 
and $b_5$ is divided by the three lines 
$x_1=x_2, \ x_1=x_3, \ x_2=x_3$ into six. 
The isotropy subgroups $W_b$ and 
the $W_b$-invariant points $z$ of $b$ 
for $b=b_1,\ldots,b_5$ 
are given in Table \ref{table:W_b-z-for-semiorder}.

\begin{table}
\caption{$W_b$ and $z$ for $b=b_1,\ldots,b_5$ 
in semiorder arrangement.}
\begin{center} 
\begin{tabular}{lll} \toprule 
$b$ & 
$W_b$ &   
$W_b$-invariant points 
$z$ 
of $b$ 
\\ \midrule 
$b_1$ & 
$\dS_{\{ 1,2\}}$ &  
$d(1,1,-2)^T, \ d>\frac{1}{3}$ 
\\ 
$b_2$ & 
$\{ 1\}$ &  
$(1,0,-1)^T+d_1(1,1,-2)^T+d_2(2,-1,-1)^T, \ 
d_1>0, \, d_2>0$ 
\\ 
$b_3$ &  
$\dS_{\{ 2,3\}}$ &  
$d(2,-1,-1)^T, \ d>\frac{1}{3}$  
\\ 
$b_4$ & 
$\{ 1\}$ &  
$(1,0,-1)^T+d_1(-2,1,1)^T+d_2(-1,-1,2)^T, \ 
d_1>0, \, d_2>0, \, d_1+d_2<\frac{1}{3}$ 
\\ 
$b_5$ & 
$\dS_{\{ 1,2,3\}}$ &  
$(0,0,0)^T$  
\\ \bottomrule 
\end{tabular}
\end{center} 
\label{table:W_b-z-for-semiorder}
\end{table} 

We can confirm 
$|W_{b_1}|=2!=|\varphi_{\calB}^{-1}(b_1)|, \ 
|W_{b_2}|=1=|\varphi_{\calB}^{-1}(b_2)|, \ 
|W_{b_3}|=2!=|\varphi_{\calB}^{-1}(b_3)|, \ 
|W_{b_4}|=1=|\varphi_{\calB}^{-1}(b_4)|, \ 
|W_{b_5}|=3!=|\varphi_{\calB}^{-1}(b_5)|$ 
(Part 2 of Theorem \ref{theorem:6}). 

We have 
$\{ W\text{-orbits of~} \Ch(\calB)\}
=\{ \calO(b_1), \ldots, \calO(b_5)\}$ and 
thus  
\[
|\{ W\text{-orbits of~} \Ch(\calB)\}|=5. 
\]  
We can see from \eqref{eq:b1b5-semiorder} that 
$|\calO(b_1)|=3, \ |\calO(b_2)|=3!, \ |\calO(b_3)|=3, \ 
|\calO(b_4)|=3\times 2=6$ and $|\calO(b_5)|=1$, 
coinciding with 
$|W|/|W_{b_i}|, \ i=1,\ldots,5$. 
With these values, 
$|\Ch(\calC)|$ can be computed also as 
\begin{eqnarray*} 
|\Ch(\calC)|
&=& \sum_{i=1}^5 |\varphi_{\calB}^{-1}(b_i)|\cdot |\calO(b_i)| 
= 2\cdot 3 +1\cdot 6+2\cdot 3+1\cdot 6+6\cdot 1 
= 30 \\ 
&=& \sum_{i=1}^5 
\left( |W_{b_i}|\times \frac{|W|}{|W_{b_i}|}\right) 
= 
|W| \times 
|\{ \text{$W$-orbits of $\Ch(\calB)$}\}|. 
\end{eqnarray*} 
In addition,  
$|\Ch(\calB)|=\sum_{i=1}^5|\calO(b_i)|
=3+6+3+6+1=19$
in agreement with \eqref{eq:DGF-|Ch(B)|}.

\begin{rem}
\label{rem:non-unique-example-semiorder}
In Table \ref{table:W_b-z-for-semiorder}, 
we find that 
$W_{b_2}=\{ 1\}$ and $W_{b_4}=\{ 1\}$ 
(resp. $W_{b_1}=\dS_{\{ 1,2\}}$ and $W_{b_3}=\dS_{\{ 2,3\}}$) 
are conjugate to each other, 
although $b_2$ and $b_4$ (resp. $b_1$ and $b_3$) 
are on different orbits.  
The map $\tau$ in \eqref{eq:def-tau} satisfies 
\[
\tau^{-1}([ \{ 1\}])
=\{ \calO(b_2), \calO(b_4)\}, 
\quad 
\tau^{-1}([ \dS_{\{ 1,2\}} ])
=\{ \calO(b_1), \calO(b_3)\}; 
\]  
hence $\tau$ is not injective.  
It is evident 
that $b_2$ and $b_4$ are on different orbits, 
because $b_2$ is a cone and $b_4$ is a triangle.  
As for $b_1$ and $b_3$, we can confirm 
$\calO(b_1)\ne \calO(b_3)$ 
by looking at their corresponding semiorders in 
\eqref{eq:b1b5-semiorder}.  
\end{rem}


\subsection{Coxeter group of type $A_{m-1}$ 
and restricted all-subset arrangement} 
\label{subsec:A+r-allsubset}

Let $W=\dS_m$ and $V=H_0$. 
Then $\calA=\calA(W)$ is the braid arrangement in $H_0$. 
%
Let $\calB$ be the restricted all-subset arrangement  
(Kamiya, Takemura and Terao \cite{ktt-aam}): 
\[
\calB=
\{ H_{I}^0 \mid 
\emptyset \ne I \subsetneq \{ 1,\ldots,m\}\}, 
\]
where 
\[
H_I^0:=\{ x=(x_1,\ldots,x_m)^T\in H_0 
\mid \sum_{i \in I}x_i=0\}. 
\]
Since $H_I^0=H_{\{ 1,\ldots,m\}\setminus I}^0$ 
for $I \ne \emptyset, \{ 1,\ldots,m\}$, 
we have $|\calB|=(2^m-2)/2=2^{m-1}-1$. 
The number of $W$-orbits of $\calB$ is 
$(m-1)/2$ if $m$ is odd and 
$m/2$ if $m$ is even. 

Theorem \ref{theorem:6} 
applied to this case gives the affirmative answer
to the open problem left in Section 6 of \cite{ktt-aam}.
Using the terminology in Corollary 6.2 of \cite{ktt-aam}, 
we state:

\begin{theorem}
\label{thm:affirmative}
The number of inequivalent ranking patterns of 
unfolding models of codimension one is 
\[
\frac{|\Ch(\calA \cup \calB)|}{m!}-1
\]
for the braid arrangement $\calA$ in $H_0$ 
and the restricted all-subset arrangement $\calB$.  
\end{theorem}

\Proof 
Part 1 of Theorem \ref{theorem:6} implies that 
the number of $W$-orbits of $\Ch(\calB)$ 
for the restricted all-subset arrangement $\calB$ 
is equal to $|\Ch(\calA \cup \calB)|/(m!).$
This completes the proof because
of the last sentence of Corollary 6.2
in \cite{ktt-aam}.
\qed 

\

In the case $m=3$, we can easily obtain 
$|\varphi_{\calA}^{-1}(a)|=2 \ (a \in \Ch(\calA)), \ 
|\Ch(\calC)|=12, \ |\Ch(\calB)|=6$ 
and other quantities by direct observations. 
Alternatively, $|\Ch(\calC)|$ can be computed by 
using the characteristic polynomial 
$\chi(\calC,t)=(t-1)(t-5)$ 
(
\cite[Sec. 6.2.1]{ktt-aam}) 
of $\calC$: 
Zaslavsky's result 
on the chamber-counting problem 
(\cite[Theorem A]{zas-75}, \cite[Theorem 2.68]{ort}) 
yields 
$|\Ch(\calC)|=(-1)^2\chi(\calC,-1)=12$.

\ 

Let us investigate the case $m=4$.  

\begin{figure}[htbp]
 \begin{center}
\includegraphics*[width=.6\textwidth]{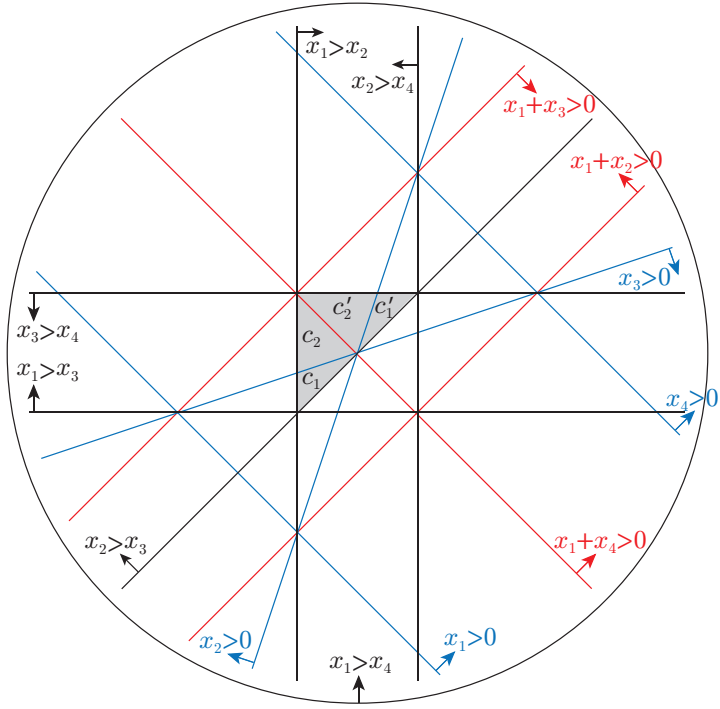}
\caption{Braid plus restricted all-subset arrangement.}  
\label{fig:braid+restricted-m=4}
 \end{center}
\end{figure}

The elements of $\calA$ are 
the 
six planes in $V=H_0, \ \dim H_0=3$, 
defined by 
the following equations:  
\begin{equation}
\label{eq:braid-m=4}
x_1=x_2, \quad x_1=x_3, \quad x_1=x_4, \quad 
x_2=x_3, \quad x_2=x_4, \quad x_3=x_4, 
\end{equation}
whereas those of $\calB$ are 
the seven planes below: 
\begin{gather}
x_1=0, \quad x_2=0, \quad x_3=0, \quad x_4=0; 
\label{eq:x1=0x4=0} \\ 
x_1+x_2=0, \quad x_1+x_3=0, \quad x_1+x_4=0. 
\label{eq:x1+x2=0x1+x4=0}
\end{gather} 
Figure \ref{fig:braid+restricted-m=4} shows the 
intersection with the unit sphere 
$\bbS^2=\{ (x_1,\ldots,x_4)^T\in H_0 \mid 
x_1^2+\cdots +x_4^2=1\}$ in $H_0$. 
Note that $\calB$ has two orbits 
under the action of 
$W$---the four planes in \eqref{eq:x1=0x4=0} 
constitute one orbit, and 
the three 
in \eqref{eq:x1+x2=0x1+x4=0}
form the other 
one.
(In Figure \ref{fig:braid+restricted-m=4}, 
the planes in \eqref{eq:x1=0x4=0} are drawn 
in blue and those in \eqref{eq:x1+x2=0x1+x4=0} are 
sketched in red.)

As a chamber of $\calA$,  
let us take $a \in \Ch(\calA)$ defined 
by $x_1>x_2>x_3>x_4$ 
($a$ is shaded in 
Figure \ref{fig:braid+restricted-m=4}).   
This chamber $a$ of $\calA$ contains 
exactly four chambers $c_1,c_2,c_1',c_2'$ 
of $\calC$: 
$\varphi_{\calA}^{-1}(a)=\{ c_1,c_2,c_1',c_2'\}$.  
These chambers have the following walls: 
\begin{eqnarray*}
&& \text{Walls of }c_1: x_1=x_2, \ x_2=x_3, \ x_3=0 
\ \ (x_1>x_2>x_3>0), \\ 
&& \text{Walls of }c_2: x_1=x_2, \ x_1+x_4=0, \ x_3=0 
\ \ (-x_4>x_1>x_2, \ x_3<0), \\ 
&& \text{Walls of }c_1': x_4=x_3, \ x_3=x_2, \ x_2=0 
\ \ (x_4<x_3<x_2<0), \\ 
&& \text{Walls of }c_2': x_4=x_3, \ x_4+x_1=0, \ x_2=0 
\ \ (-x_1<x_4<x_3, \ x_2>0). 
\end{eqnarray*} 
Note that $c_1'$ (resp. $c_2'$) is obtained 
from $c_1$ (resp. $c_2$) 
by changing $(x_1,x_2,x_3,x_4)$ to 
$(-x_4,-x_3,-x_2,-x_1)$. 
Since 
$|\varphi_{\calA}^{-1}(a)|=4$, 
we have $|\Ch(\calC)|
=|\varphi_{\calA}^{-1}(a)|\times |W|=4\times 4!=96$, 
and $\Ch(\calC)$ has exactly four orbits 
$\calO(c_1), \calO(c_2), \calO(c_1'), \calO(c_2')$. 

\begin{figure}[htbp]
 \begin{center}
\includegraphics*[width=.6\textwidth]{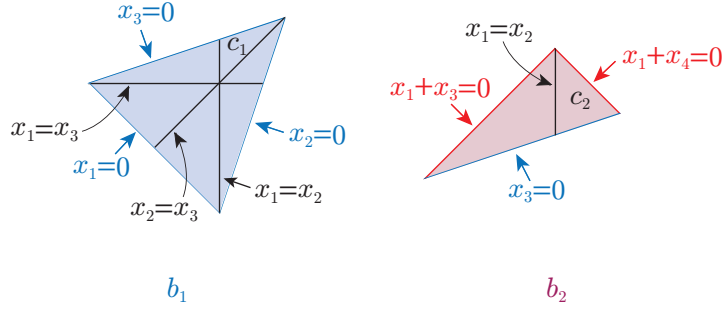}
\caption{Braid plus restricted all-subset arrangement. 
$b_1$ and $b_2$.}
\label{fig:braid+restricted-m=4-b}
 \end{center}
\end{figure}

The chambers 
$b_i:=\varphi_{\calB}(c_i), \ 
b_i':=\varphi_{\calB}(c_i'), \ i=1,2$, 
of $\calB$ containing $c_1,c_2,c_1',c_2'$ 
have the following walls: 
\begin{eqnarray*}
&& \text{Walls of }b_1: x_1=0, \ x_2=0, \ x_3=0 
\ \ (x_1>0, \ x_2>0, \ x_3>0), \\ 
&& \text{Walls of }b_2: x_1+x_3=0, \ x_1+x_4=0, \ x_3=0 
\ \ (-x_4>x_1>-x_3>0), \\ 
&& \text{Walls of }b_1': x_4=0, \ x_3=0, \ x_2=0 
\ \ (x_4<0, \ x_3<0, \ x_2<0), \\ 
&& \text{Walls of }b_2': x_4+x_2=0, \ x_4+x_1=0, \ x_2=0 
\ \ (-x_1<x_4<-x_2<0) 
\end{eqnarray*} 
(Figure \ref{fig:braid+restricted-m=4-b}). 
The chamber $b_1\in \Ch(\calB)$ is divided 
by the three planes 
$x_1=x_2, \ x_1=x_3, \ x_2=x_3$ 
of $\calA$ 
into six chambers of $\calC$, 
whereas $b_2$ is divided by the plane 
$x_1=x_2$ into two chambers of $\calC$.  
For $b_1$, we have 
$W_{b_1}=\dS_{\{ 1,2,3\}}$ 
(the $W_{b_1}$-invariant points $z$ of $b_1$ are 
$z=d(1,1,1,-3)^T, \ d>0$), 
and for $b_2$, we find $W_{b_2}=\dS_{\{ 1,2\}}$  
(the $W_{b_2}$-invariant points $z$ of $b_2$ are 
$z=d_1(1,1,-1,-1)^T+d_2(1,1,0,-2)^T, \ d_1,d_2>0$). 
So we see 
$|W_{b_1}|=|\varphi_{\calB}^{-1}(b_1)| \ (=6)$ and 
$|W_{b_2}|=|\varphi_{\calB}^{-1}(b_2)| \ (=2)$ 
hold true.

We have 
$\{ W\text{-orbits of~} \Ch(\calB)\}
=\{ \calO(b_1), \calO(b_2), \calO(b_1'),\calO(b_2')\}$ 
and hence
\[
|\{ W\text{-orbits of~} \Ch(\calB)\}|=4. 
\] 
The chambers $b \in \Ch(\calB)$ on the same $W$-orbit 
as $b_1$, $\calO(b)=\calO(b_1)$, have walls of the form 
$x_i>0, \ x_j>0, \ x_k>0$ 
($i,j,k$ are all distinct),  
while $b\in \Ch(\calB)$ such that 
$\calO(b)=\calO(b_2)$ have walls  
of the form 
$x_i<0, \ x_i+x_j>0, \ x_i+x_k>0$ 
($i,j,k$ are all distinct).    
Thus the orbit sizes are 
$|\calO(b_1)|=\binom{4}{3}=4 \ (=|W|/|W_{b_1}|)$ and 
$|\calO(b_2)|=\binom{4}{3}\times 3=12  
\ (=|W|/|W_{b_2}|)$.   
Accordingly, $|\Ch(\calC)|$ is again 
\[
|\Ch(\calC)|
=\left(|\varphi_{\calB}^{-1}(b_1)|\cdot |\calO(b_1)|
+|\varphi_{\calB}^{-1}(b_2)|\cdot |\calO(b_2)|\right)\times 2 
= (6\times 4+2\times 12)\times 2 = 96. 
\] 
Besides, 
$|\Ch(\calB)|=2(|\calO(b_1)|+|\calO(b_2)|)
=2(4+12)=32$.  

The characteristic polynomial 
$\chi(\calC,t)$ of $\calC$ 
is  
\begin{eqnarray*}
\chi(\calC,t)=(t-1)(t-5)(t-7)
\end{eqnarray*} 
(Kamiya,Takemura and Terao 
\cite[Sec. 6.2.2]{ktt-aam}). 
This polynomial yields 
\[
|\Ch(\calC)|=(-1)^3\chi(\calC,-1)=96 
\] 
in agreement with our observations above.   

\ 

For $5 \le m\le 9$, we used the finite-field method 
(Athanasiadis \cite{ath-thesis, ath-96}, 
Stanley \cite[Lecture 5]{sta-07}, 
Crapo and Rota \cite{crr}, 
Kamiya, Takemura and Terao \cite{ktt-08, ktt-root, 
ktt-nc}) 
to calculate the characteristic polynomials 
$\chi(\calC,t)$ of $\calC$, 
and obtained the numbers of $W$-orbits of $\Ch(\calB)$ 
by using Zaslavsky's result \cite[Theorem A]{zas-75} 
and Part 1 of Theorem \ref{theorem:6}: 
$|\{ \text{$W$-orbits of $\Ch(\calB)$}\}|
=|\varphi_{\calA}^{-1}(a)|=|\Ch(\calC)|/|W|$ as follows.    
\begin{align*}
m=5: \quad \chi(\calC,t)&=(t-1)(t-7)(t-8)(t-9), 
\quad |\Ch(\calC)|=1440, \quad |\varphi_{\calA}^{-1}(a)|=12, \\  
m=6: \quad \chi(\calC,t)&=(t-1)(t-7)(t-11)(t-13)(t-14), \\ 
& |\Ch(\calC)|=40320, \quad |\varphi_{\calA}^{-1}(a)|=56, \\  
m=7: \quad 
\chi(\calC,t)&=(t-1)(t-11)(t-13)(t-17)(t-19)(t-23), \\
& |\Ch(\calC)|=2903040, \quad |\varphi_{\calA}^{-1}(a)|=576, \\  
m=8: \quad \chi(\calC,t)&= (t-1)(t-19)(t-23)(t-25)(t-27)(t-29)(t-31), \\ 
& |\Ch(\calC)|=670924800, \quad |\varphi_{\calA}^{-1}(a)|=16640, \\  
m=9: \quad \chi(\calC,t)&=
(t-1)(t^7-290t^6+36456t^5-2573760t^4+110142669t^3 \\ 
& \qquad \qquad \qquad -2855339970t^2+41492561354t-260558129500), \\ 
& |\Ch(\calC)|=610037568000, \quad |\varphi_{\calA}^{-1}(a)|=1681100.
\end{align*}

Note that the characteristic polynomial $\chi(\calC, t)$ factors 
into polynomials of degree one over $\bbZ$ for $m \le 8$.  

\begin{rem}
For $m=5$ and $m=6$, 
Kamiya, Takemura and Terao \cite{ktt-aam} 
identified all the elements $c$ of 
$\varphi_{\calA}^{-1}(a)$ for 
$a:x_1>\cdots >x_m$ 
and gave 
an example 
of 
the $W_b$-invariant 
points 
$z$ of $b=\varphi_{\calB}(c)$ for each $c$. 
From those $z$, we immediately obtain $W_b$ 
by Proposition \ref{prop:W_b}. 
\end{rem}

\subsection{Coxeter group of type $A_{m-1}$ 
and unrestricted all-subset arrangement} 
\label{subsec:A+u-allsubset}

Let $W=\dS_m$ and $V=\bbR^m$. 
Then $\calA=\calA(W)$ is the braid arrangement 
in $\bbR^m$. 
Let $\calB$ be the 
(unrestricted) all-subset arrangement 
(\cite{ktt-aam}): 
\[
\calB=\{ H_{I} \mid 
\emptyset \ne I \subseteq \{ 1,\ldots,m\}\}, 
\]
where 
\[
H_I:=\{ x=(x_1,\ldots,x_m)^T\in \bbR^m 
\mid \sum_{i \in I}x_i=0\}. 
\]
Note $|\calB|=2^m-1$. 
The number of orbits of $\calB$ under 
the action of $W$ is $m$. 

\ 

We will examine the  case $m=3$.  

\begin{figure}[htbp]
 \begin{center}
\includegraphics*[width=.55\textwidth]{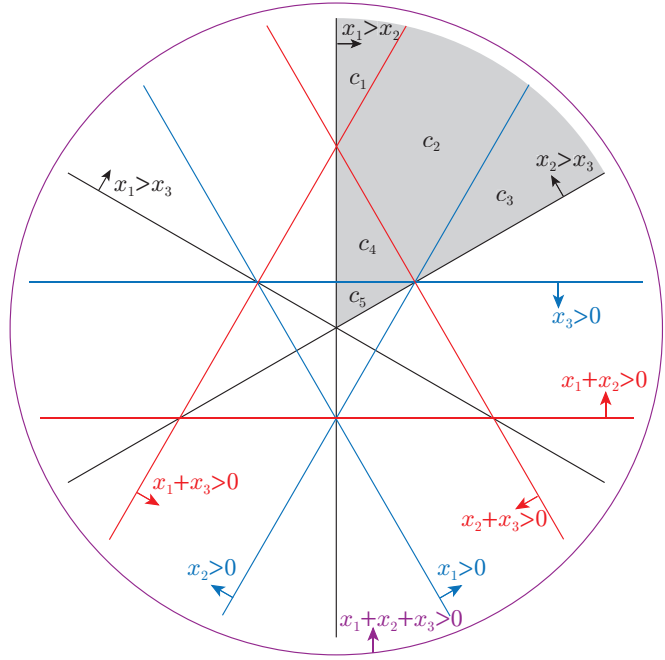}
\caption{Braid plus unrestricted all-subset arrangement.}
\label{fig:braid+unrestricted-m=3}
 \end{center}
\end{figure}

The arrangement $\calA$ has 
exactly 
the three planes in 
$V=\bbR^3$ defined by the same equations as 
those in \eqref{eq:braid-with-m=3}. 
On the other hand, $\calB$ consists of the seven planes 
defined by 
\begin{gather*}
\label{eq:u-all-m=3}
x_1=0, \quad x_2=0, \quad x_3=0; \\ 
x_1+x_2=0, \quad x_1+x_3=0, \quad x_2+x_3=0; \\ 
x_1+x_2+x_3=0
\end{gather*}
with each line corresponding to one orbit 
under the action of $W$ on $\calB$. 
Figure \ref{fig:braid+unrestricted-m=3} exhibits 
the intersection with 
the unit sphere in $V=\bbR^3$. 

Let us take $a \in \Ch(\calA)$ defined 
by $x_1>x_2>x_3$ 
($a$ with $x_1+x_2+x_3>0$ is shaded in 
Figure \ref{fig:braid+unrestricted-m=3}). 
Then 
$\allowbreak \varphi_{\calA}^{-1}(a)
=\{ c_1,c_2,c_3,c_4,c_5,c_1',c_2',c_3',c_4',c_5'\}$,  
where 
\begin{eqnarray*}
&& c_1: x_1>x_2, \ x_1+x_3<0, \ x_1+x_2+x_3>0, \\ 
&& c_2: x_2>0, \ x_1+x_3>0, \ x_2+x_3<0, \\ 
&& c_3: x_2>x_3, \ x_2<0, \ x_1+x_2+x_3>0, \\ 
&& c_4: x_1>x_2, \ x_3<0, \ x_2+x_3>0, \\ 
&& c_5: x_1>x_2, \ x_2>x_3, \ x_3>0, 
\end{eqnarray*} 
and $c_i', \ i=1,\ldots,5$, are the chambers 
obtained from 
$c_i, \ i=1,\ldots,5$, by changing $(x_1,x_2,x_3)$ to 
$(-x_3,-x_2,-x_1)$.  
Thus 
$|\varphi_{\calA}^{-1}(a)|
=|\{ \text{$W$-orbits of $\Ch(\calC)$}\}|
=10$ 
and 
$|\Ch(\calC)|=|\varphi_{\calA}^{-1}(a)|\times |W|
=10\times 3!=60$.   

\begin{figure}[htbp]
 \begin{center}
\includegraphics*[width=.65\textwidth]{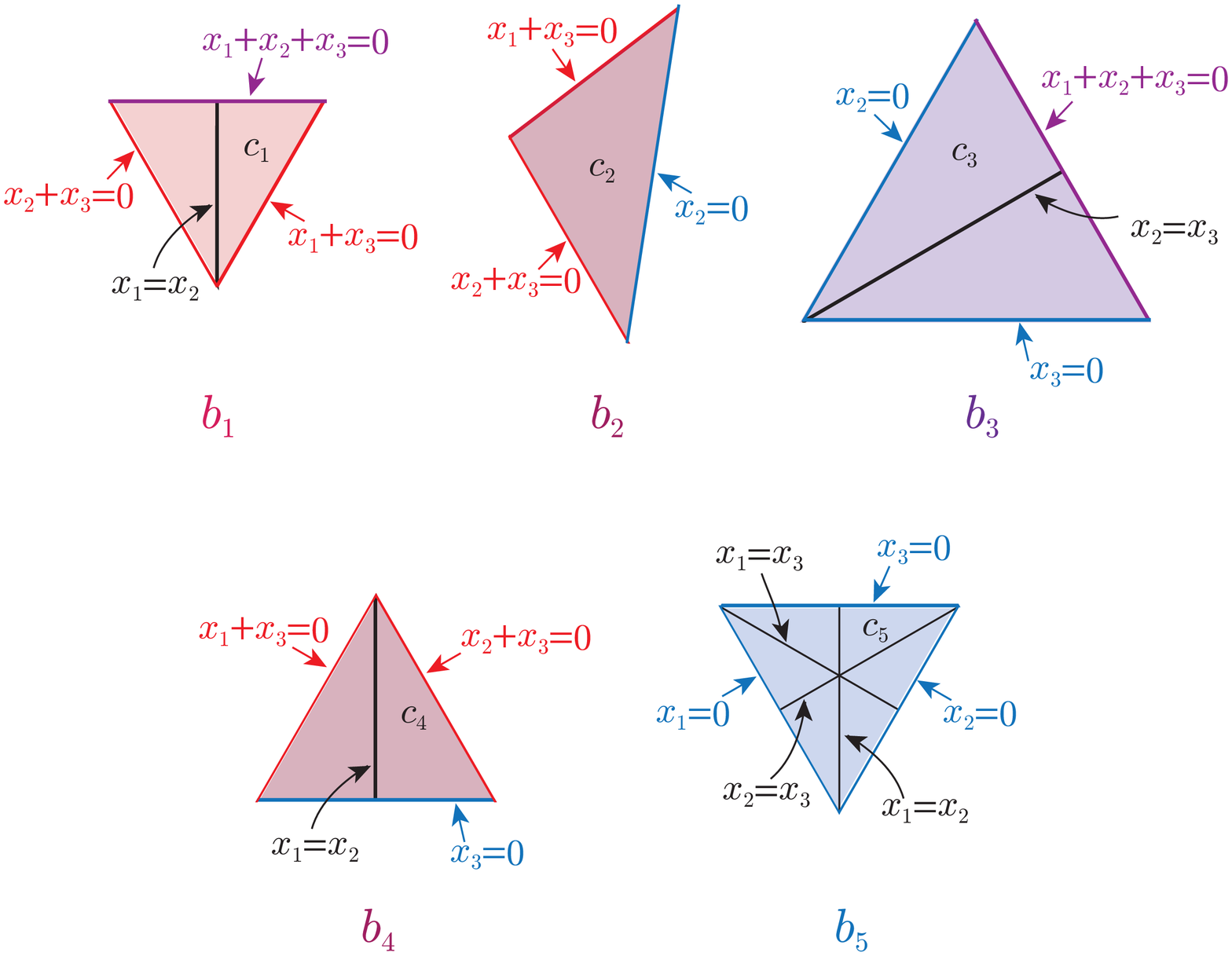}
\caption{Braid plus unrestricted all-subset arrangement. $b_1,\ldots,b_5$.}
\label{fig:braid+unrestricted-m=3-b}
 \end{center}
\end{figure}

The chambers 
$b_i:=\varphi_{\calB}(c_i) \in \Ch(\calB), 
\ i=1,\ldots,5$, 
are 
\begin{eqnarray}
\label{eq:b1-b5}
&& b_1: x_1+x_3<0, \ x_2+x_3<0, \ x_1+x_2+x_3>0, 
\nonumber \\ 
&& b_2: x_2>0, \ x_1+x_3>0, \ x_2+x_3<0, \nonumber \\ 
&& b_3: x_2<0, \ x_3<0, \ x_1+x_2+x_3>0, \\ 
&& b_4: x_3<0, \ x_1+x_3>0, \ x_2+x_3>0, \nonumber \\ 
&& b_5: x_1>0, \ x_2>0, \ x_3>0, \nonumber 
\end{eqnarray} 
and $b_i':=\varphi_{\calB}(c_i') \in \Ch(\calB), 
\ i=1,\ldots,5$, 
can be obtained from $b_i, \ i=1,\ldots,5$,  
by the above-mentioned rule. 
See Figure \ref{fig:braid+unrestricted-m=3-b}. 
The chamber $b_1$ is divided 
by the plane $x_1=x_2$ into two chambers; 
$b_2$ is not divided by any 
plane 
in $\calA$; 
$b_3$ is divided by $x_2=x_3$ into two;  
$b_4$ is divided by $x_1=x_2$ into two;  
and $b_5$ is divided by the three planes 
$x_1=x_2, \ x_1=x_3, \ x_2=x_3$ into six.  
The isotropy subgroups $W_b$ and 
the $W_b$-invariant points $z$ of $b$ 
for $b=b_1,\ldots,b_5$ 
are given in Table \ref{table:W_b-z-for-u-all-subset}. 
\begin{table}
\caption{$W_b$ and $z$ for $b=b_1,\ldots,b_5$ 
in unrestricted all-subset arrangement.}
\begin{center} 
\begin{tabular}{lll} \toprule 
$b$ & 
$W_b$ &   
$W_b$-invariant points 
$z$ 
of $b$ 
\\ \midrule 
$b_1$ & 
$\dS_{\{ 1,2\}}$ &  
$d_1(1,1,-1)^T+d_2(1,1,-2)^T, \ d_1,d_2>0$ 
\\ 
$b_2$ & 
$\{ 1\}$ &  
$d_1(1,1,-1)^T+d_2(1,0,-1)^T+d_3(1,0,0)^T, \ 
d_1,d_2,d_3>0$ 
\\ 
$b_3$ & 
$\dS_{\{ 2,3\}}$ &  
$d_1(1,0,0)^T+d_2(2,-1,-1)^T, \ d_1,d_2>0$  
\\ 
$b_4$ & 
$\dS_{\{ 1,2\}}$ &  
$d_1(1,1,-1)^T+d_2(1,1,0)^T, \ d_1,d_2>0$ 
\\ 
$b_5$ &  
$\dS_{\{ 1,2,3\}}$ & 
$d(1,1,1)^T, \ d>0$  
\\ \bottomrule 
\end{tabular}
\end{center} 
\label{table:W_b-z-for-u-all-subset}
\end{table}
We can confirm 
$|W_{b_1}|=2!=|\varphi_{\calB}^{-1}(b_1)|, \ 
|W_{b_2}|=1=|\varphi_{\calB}^{-1}(b_2)|, \ 
|W_{b_3}|=2!=|\varphi_{\calB}^{-1}(b_3)|, \ 
|W_{b_4}|=2!=|\varphi_{\calB}^{-1}(b_4)|, \ 
|W_{b_5}|=3!=|\varphi_{\calB}^{-1}(b_5)|$.

We have 
$\{ W\text{-orbits of~} \Ch(\calB)\}
=\{ \calO(b_1), \ldots, \calO(b_5),
\calO(b_1'), \ldots, \calO(b_5')\}$ and 
thus  
\[
|\{ W\text{-orbits of~} \Ch(\calB)\}|=10. 
\] 
From \eqref{eq:b1-b5}, we see 
$|\calO(b_1)|=3, \ |\calO(b_2)|=3\times 2=6, \ 
|\calO(b_3)|=3, \ |\calO(b_4)|=3, \ |\calO(b_5)|=1$, 
which coincide with 
$|W|/|W_{b_i}|, \ i=1,\ldots,5$. 
Hence, $|\Ch(\calC)|=60$ can be obtained also from 
\[
|\Ch(\calC)|
=2\sum_{i=1}^5
|\varphi_{\calB}^{-1}(b_i)|\cdot |\calO(b_i)| 
= 2 \, (2\cdot 3+1\cdot 6+2\cdot 3+2\cdot 3
+6\cdot 1) 
= 60. 
\]
We can also get 
$|\Ch(\calB)|=2\sum_{i=1}^5|\calO(b_i)|
=2\,(3+6+3+3+1)=32$.  


\begin{rem}
\label{rem:non-unique-example}
In Table \ref{table:W_b-z-for-u-all-subset}, 
we find that 
$W_{b_1}=W_{b_4}=\dS_{\{1,2\}}, \ 
W_{b_3}=\dS_{\{2,3\}}$ are all conjugate to one another, 
although $b_1, b_3, b_4$ are on different orbits.  
We have 
\[
\tau^{-1}([ \dS_{\{1,2\}}])
=\{ \calO(b_1), \calO(b_3), \calO(b_4)\}, 
\]  
so $\tau$ is not injective.  
The chambers $b_1, b_3, b_4$ are 
triangular cones 
(triangles in Figure \ref{fig:braid+unrestricted-m=3-b}) 
cut by a single plane (line) $x_i=x_j$ 
from the braid arrangement.  
However, these chambers are 
easily seen to be 
on different orbits, 
since their three walls (edges)  
are of different combinations of orbits of $\calB$.
\end{rem}

For $m\le 7$, we computed $\chi(\calC,t)$ using 
the finite-field method, 
and obtained the numbers of $W$-orbits of $\Ch(\calB)$ 
as follows:   
\begin{align*}
m=3: \quad \chi(\calC,t) &=(t-1)(t-4)(t-5), 
\quad |\Ch(\calC)|=60, \quad 
|\varphi_{\calA}^{-1}(a)|=10, \\  
m=4: \quad \chi(\calC,t) &=(t-1)(t-5)(t-7)(t-8), 
\quad |\Ch(\calC)|=864, \quad 
|\varphi_{\calA}^{-1}(a)|=36, \\  
m=5: \quad \chi(\calC,t) &=(t-1)(t-7)(t-9)(t-11)(t-13), \\ 
& |\Ch(\calC)|=26880, \quad 
|\varphi_{\calA}^{-1}(a)|=224, \\  
m=6: \quad \chi(\calC,t) &=(t-1)(t-11)(t-13)(t-17)^2(t-19), \\ 
& |\Ch(\calC)|=2177280, \quad 
|\varphi_{\calA}^{-1}(a)|=3024, \\  
m=7: \quad \chi(\calC,t) &=(t-1)(t-19)(t-23)(t^4-105 t^3+4190 t^2-75180 t + 510834), \\
& |\Ch(\calC)|=566697600, \quad 
|\varphi_{\calA}^{-1}(a)|=112440. 
\end{align*}

Note that the characteristic polynomial $\chi(\calC, t)$ factors 
into polynomials of degree one over $\bbZ$ for $m \le 6$.

\subsection{Mid-hyperplane arrangement} 
\label{subsec:mid-hyper}

Let $W=\dS_m$ and 
$V=H_0$, 
so $\calA=\calA(W)$ is the braid arrangement 
in $H_0$. 
We take  
\[
\calB=\{ H_{ijkl} 
\mid 1\le i<j\le m, \ 1\le k<l\le m, \ i<k, \ 
|\{ i,j,k,l\}|=4 \}, 
\]
where 
\[
H_{ijkl}:=
\{ x=(x_1,\ldots,x_m)^T \in 
H_0 
\mid  x_i+x_j=x_k+x_l\}, 
\]
so that $\calC=\calA\cup \calB$ is 
an essentialization of 
the mid-hyperplane arrangement 
(Kamiya, Orlik, Takemura and Terao \cite{kott}).  
We have $|\calB|=3\binom{m}{4}$,  
and the action of $W$ on $\calB$ is transitive.  

\ 

Let us consider the case $m=4$.  

\begin{figure}[htbp]
 \begin{center}
\includegraphics*[width=.5\textwidth]{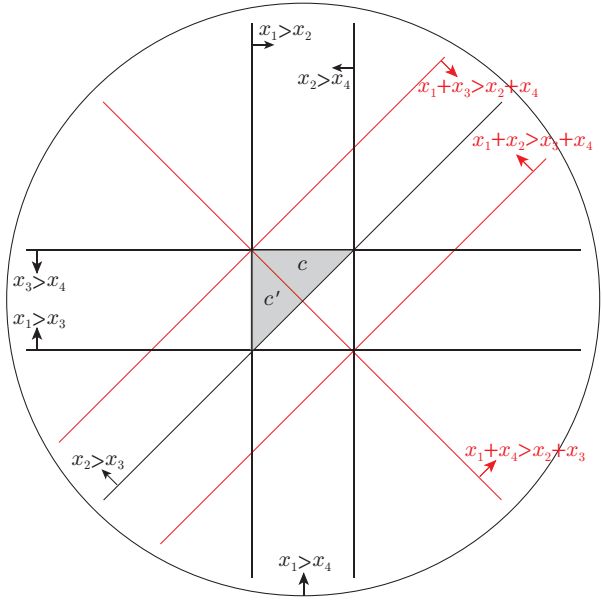}
\caption{Essentialization of 
mid-hyperplane arrangement.}
\label{fig:mid-hyper-m=4}
 \end{center}
\end{figure}

The elements of $\calA$ are the six planes 
in $V=H_0, \ \dim H_0=3$, defined by the 
equations in \eqref{eq:braid-m=4}, 
whereas those of $\calB$ are the three planes 
defined by the following equations:  
\[
x_1+x_2=x_3+x_4, \quad x_1+x_3=x_2+x_4, \quad 
x_1+x_4=x_2+x_3. 
\]
Figure \ref{fig:mid-hyper-m=4} shows 
the intersection with 
the unit sphere $\bbS^2$ in $H_0$. 

Let us take $a \in \Ch(\calA)$ defined 
by $x_1>x_2>x_3>x_4$ 
($a$ is shaded in Figure \ref{fig:mid-hyper-m=4}). 
Then 
$\allowbreak \varphi_{\calA}^{-1}(a)
=\{ c,c'\}$,  
where 
\begin{eqnarray*}
&& c: x_2>x_3, \ x_3>x_4, \ x_1+x_4>x_2+x_3, \\ 
&& c': x_3<x_2, \ x_2<x_1, \ x_4+x_1<x_3+x_2.  
\end{eqnarray*}  
Note that $c'$ is obtained from $c$ by changing 
$(x_1,x_2,x_3,x_4)$ to $(-x_4,-x_3,-x_2,-x_1)$. 
We have   
$|\varphi_{\calA}^{-1}(a)|
=|\{ \text{$W$-orbits of $\Ch(\calC)$}\}|
=2$ 
and 
$|\Ch(\calC)|=|\varphi_{\calA}^{-1}(a)|\times |W|
=2\times 4!=48$.   

\begin{figure}[htbp]
 \begin{center}
\includegraphics*[width=.35\textwidth]{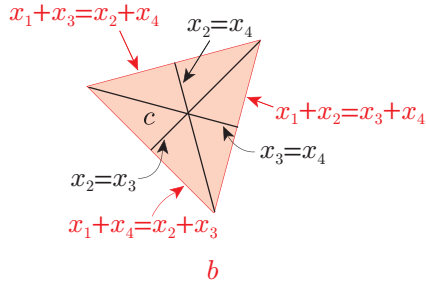}
\caption{Essentialization of 
mid-hyperplane arrangement. $b$.}
\label{fig:mid-hyper-m=4-b}
 \end{center}
\end{figure}

The chamber 
$b:=\varphi_{\calB}(c) \in \Ch(\calB)$ 
is 
\begin{equation}
\label{eq:b:x1+x2>x3+x4}
b: 
x_1+x_2>x_3+x_4, \ 
x_1+x_3>x_2+x_4, \ 
x_1+x_4>x_2+x_3,  
\end{equation}
which is divided 
by the three planes 
$x_2=x_3, \ x_2=x_4, \ x_3=x_4$ into six chambers  
(Figure \ref{fig:mid-hyper-m=4-b}).  
We find 
the $W_b$-invariant points $z$ of $b$ to be 
$z=d(3,-1,-1,-1)^T, \ d>0$, 
so 
$W_{b}=\dS_{\{ 2,3,4\}}$ 
and 
$|W_{b}|=3!=|\varphi_{\calB}^{-1}(b)|$.

We have 
$\{ W\text{-orbits of~} \Ch(\calB)\}
=\{ \calO(b), \calO(b')\}$ 
with $b':=\varphi_{\calB}(c')$. 
Thus  
\[
|\{ W\text{-orbits of~} \Ch(\calB)\}|=2. 
\] 
From \eqref{eq:b:x1+x2>x3+x4}, we see 
$|\calO(b)|=4 \ 
(=|W|/|W_{b}|=4!/(3!))$, 
so we can calculate $|\Ch(\calC)|$ alternatively as 
\[
|\Ch(\calC)|
=2(|\varphi_{\calB}^{-1}(b)|\cdot |\calO(b)|) 
= 2\times 6\cdot 4 
= 48. 
\]
Moreover, we get 
$|\Ch(\calB)|=2|\calO(b)|=2\times 4=8$.  


\ 

For $m\le 10$, the characteristic polynomials 
of $\calC$ are known 
(\cite{kott}, \cite{ktts-appl}), 
so we can find the numbers of $W$-orbits of $\Ch(\calB)$:   
\begin{align*}
m=4: \quad \chi(\calC,t) &= 
(t-1)(t-3)(t-5), \quad 
|\Ch(\calC)|=48, \quad 
|\varphi_{\calA}^{-1}(a)|=2, \\  
m=5: \quad \chi(\calC,t) &= 
(t-1)(t-7)(t-8)(t-9), \quad 
|\Ch(\calC)|=1440, \quad 
|\varphi_{\calA}^{-1}(a)|=12, \\  
m=6: \quad \chi(\calC,t) &= 
(t-1)(t-13)(t-14)(t-15)(t-17), \\ 
& |\Ch(\calC)|=120960, \quad 
|\varphi_{\calA}^{-1}(a)|=168, \\  
m=7: \quad \chi(\calC,t) &= 
(t-1)(t-23)(t-24)(t-25)(t-26)(t-27), \\ 
& |\Ch(\calC)|=23587200, \quad 
|\varphi_{\calA}^{-1}(a)|=4680, \\  
m=8: \quad \chi(\calC,t) &= 
(t-1)(t-35)(t-37)(t-39)(t-41)(t^2-85t+1926), \\ 
& |\Ch(\calC)| = 9248117760, \quad 
|\varphi_{\calA}^{-1}(a)| = 229386, \\ 
m=9: \quad \chi(\calC, t) &= 
(t-1)(t^7 - 413t^6 + 73780t^5 - 7387310t^4 
+ 447514669t^3 \\ 
& \qquad \qquad \qquad 
- 16393719797t^2 + 336081719070t - 2972902161600), \\ 
& |\Ch(\calC)|=6651665153280, \quad 
|\varphi_{\calA}^{-1}(a)| = 18330206, \\ 
m=10: \quad \chi(\calC,t) &= 
(t - 1)(t^8 - 674t^7 + 201481t^6 - 34896134t^5 
+ 3830348179t^4 \\ 
& \qquad \qquad \qquad 
- 272839984046t^3 + 12315189583899t^2 \\ 
& \qquad \qquad \qquad \qquad 
- 321989533359786t + 3732690616086600), \\ 
& |\Ch(\calC)|=
8134544088921600, \quad 
|\varphi_{\calA}^{-1}(a)| = 
2241662282. 
\end{align*}

\subsection{
Signed 
all-subset arrangement} 
\label{subsec:signed-allsubset}

Let $W$ be the Coxeter group of type $B_m$, 
i.e., the semidirect product of $\dS_m$ by 
$(\bbZ/2\bbZ)^m$: $W=(\bbZ/2\bbZ)^m\rtimes \dS_m, 
\ |W|=2^m \cdot m!$. 
Then $W$ acts on $V=\bbR^m$ 
by permuting coordinates by $\dS_m$ 
and changing signs of coordinates by $(\bbZ/2\bbZ)^m$. 
The Coxeter arrangement $\calA=\calA(W)$ 
consists of the hyperplanes defined by 
\begin{gather}
x_i=0, \quad 1\le i\le m;  
\label{eq:xi=0} \\ 
x_i+x_j=0, \ x_i-x_j=0, \quad 1\le i<j\le m. 
\label{eq:xi+xj=0xi-xj=0} 
\end{gather}
We have 
$|\calA|=m+m(m-1)=m^2$. 
Moreover, the number of orbits of $\calA$ 
under the action of $W$ is two: 
one consisting of the $m$ hyperplanes in 
\eqref{eq:xi=0} and the other made up of 
the $m(m-1)$ 
hyperplanes 
in \eqref{eq:xi+xj=0xi-xj=0}. 

Let 
\[
\calB=\{ H_{(\epsilon_1,\ldots,\epsilon_m)} 
\mid \epsilon_1,\ldots,\epsilon_m \in \{ -1,0,1\}, 
\ \sum_{i=1}^m|\epsilon_i| \ge 3\}, 
\]
where 
\[
H_{(\epsilon_1,\ldots,\epsilon_m)}:=
\{ x=(x_1,\ldots,x_m)^T \in \bbR^m \mid  
\sum_{i=1}^m\epsilon_i x_i=0\}. 
\]
Note $H_{(\epsilon_1,\ldots,\epsilon_m)}=
H_{(-\epsilon_1,\ldots,-\epsilon_m)}$ 
so that 
$|\calB|=\sum_{i=3}^{m}2^{i-1}\binom{m}{i}$.  
The number of $W$-orbits of $\calB$ is $m-2$.

\ 

Let us study the case $m=3$.  

\begin{figure}[htbp]
 \begin{center}
\includegraphics*[width=.55\textwidth]{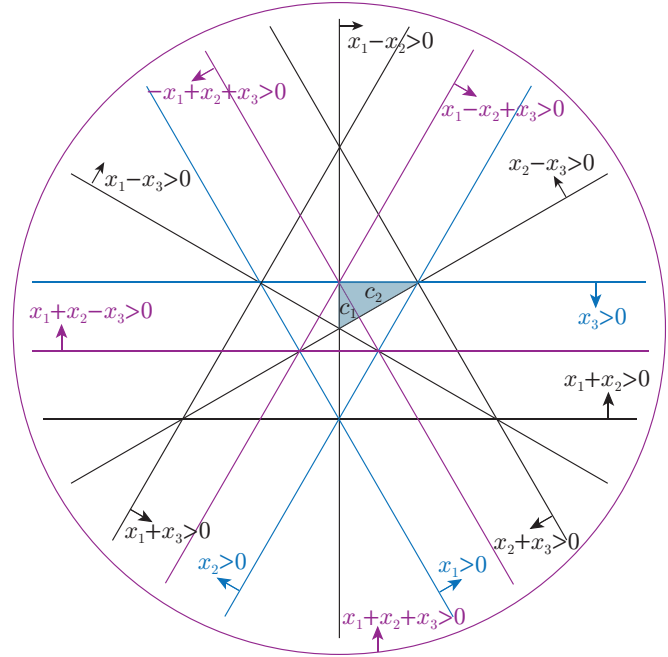}
\caption{Signed all-subset arrangement.}
\label{fig:signed-m=3}
 \end{center}
\end{figure}

In this case, $\calA$ comprises the nine 
planes in $V=\bbR^3$ defined by 
\begin{gather}
x_1=0, \ x_2=0, \ x_3=0; 
\label{eq:x1=0x3=0} \\ 
x_1+x_2=0, \ x_1-x_2=0, \ x_1+x_3=0, \ x_1-x_3=0, \ 
x_2+x_3=0, \ x_2-x_3=0 
\label{eq:x1+x2=0x2-x3=0}
\end{gather}
with each line corresponding to one orbit, 
and $\calB$ consists of one orbit 
containing 
the four planes defined by  
\begin{equation}
\label{eq:-x1+x2+x3=0}
-x_1+x_2+x_3=0, \ x_1-x_2+x_3=0, \ x_1+x_2-x_3=0, \ 
x_1+x_2+x_3=0. 
\end{equation}
Figure \ref{fig:signed-m=3} shows the intersection with 
the unit sphere in $V=\bbR^3$. 
(In Figure \ref{fig:signed-m=3}, 
the planes in \eqref{eq:x1=0x3=0}, 
\eqref{eq:x1+x2=0x2-x3=0} and \eqref{eq:-x1+x2+x3=0} 
are drawn in blue, black and purple, respectively.)

Let us take $a \in \Ch(\calA)$ defined 
by 
$x_1>x_2>x_3>0$  
(this chamber is shaded in 
Figure \ref{fig:signed-m=3}). 
Then 
$\allowbreak \varphi_{\calA}^{-1}(a)
=\{ c_1,c_2\}$,  
where 
\begin{eqnarray*}
&& c_1: x_1-x_2>0, \ x_2-x_3>0, \ -x_1+x_2+x_3>0, \\ 
&& c_2: x_3>0, \ x_2-x_3>0, \ -x_1+x_2+x_3<0.  
\end{eqnarray*} 
So  
$|\varphi_{\calA}^{-1}(a)|
=|\{ \text{$W$-orbits of $\Ch(\calC)$}\}|
=2$ 
and 
$|\Ch(\calC)|=|\varphi_{\calA}^{-1}(a)|\times |W|
=2\times 2^3\cdot 3!=96$.   

\begin{figure}[htbp]
 \begin{center}
\includegraphics*[width=.6\textwidth]{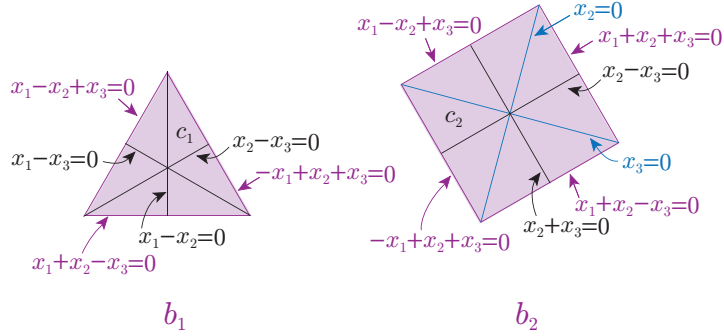}
\caption{Signed all-subset arrangement. 
$b_1$ and $b_2$.}
\label{fig:signed-m=3-b}
 \end{center}
\end{figure}

The chambers 
$b_i:=\varphi_{\calB}(c_i) \in \Ch(\calB), 
\ i=1,2$, 
are 
\begin{eqnarray*}
&& b_1: -x_1+x_2+x_3>0, \ x_1-x_2+x_3>0, \ 
x_1+x_2-x_3>0, \\ 
&& b_2: -x_1+x_2+x_3<0, \ x_1-x_2+x_3>0, \ 
x_1+x_2-x_3>0, \ x_1+x_2+x_3>0  
\end{eqnarray*}  
(Figure \ref{fig:signed-m=3-b}). 
The chamber $b_1$ is divided 
by the three planes $x_1-x_2=0, \ x_1-x_3=0, \ 
x_2-x_3=0$ into six chambers, 
and $b_2$ is divided by the four planes 
$x_2=0, \ x_3=0, \ x_2+x_3=0, \ x_2-x_3=0$ 
into eight.  
We see 
$W_{b_1}=\dS_{\{ 1,2,3\}}$ 
is the Coxeter group of type $A_2$ 
(the $W_{b_1}$-invariant points $z$ of $b_1$ are 
$z=d(1,1,1)^T, \ d>0$), 
and that 
$W_{b_2}=(\bbZ/2\bbZ)^2\rtimes \dS_{\{ 2,3\}}$ 
is the Coxeter group of type $B_2$ 
(the $W_{b_2}$-invariant points $z$ of $b_2$ are 
$z=d(1,0,0)^T, \ d>0$). 
Hence  
$|W_{b_1}|=3!=|\varphi_{\calB}^{-1}(b_1)|, \ 
|W_{b_2}|=2^2\cdot 2!=|\varphi_{\calB}^{-1}(b_2)|$.

We have 
$\{ W\text{-orbits of~} \Ch(\calB)\}
=\{ \calO(b_1), \calO(b_2)\}$,  
so  
\[
|\{ W\text{-orbits of~} \Ch(\calB)\}|=2. 
\] 

From Figures \ref{fig:signed-m=3} and 
\ref{fig:signed-m=3-b}, we see 
$|\calO(b_1)|=4\times 2=8 \ 
(=|W|/|W_{b_1}|=(2^3\cdot 3!)/(3!)), 
\ |\calO(b_2)|=3\times 2=6 \ 
(=|W|/|W_{b_2}|=(2^3\cdot 3!)/(2^2\cdot 2!))$, 
so $|\Ch(\calC)|$ can be computed also as 
\[
|\Ch(\calC)|
=|\varphi_{\calB}^{-1}(b_1)|\cdot |\calO(b_1)|
+ |\varphi_{\calB}^{-1}(b_2)|\cdot |\calO(b_2)|
= 6\cdot 8+8\cdot 6
= 96. 
\]
Furthermore, we can get 
$|\Ch(\calB)|=|\calO(b_1)|+|\calO(b_2)|=8+6=14$.  



\ 

For $m \le 6$, we computed $\chi(\calC,t)$ and 
obtained the numbers of $W$-orbits of $\Ch(\calB)$ 
as follows:  
\begin{align*}
m=3: \quad 
\chi(\calC,t)&=(t-1)(t-5)(t-7), \quad 
|\Ch(\calC)| = 96, \quad 
|\varphi_{\calA}^{-1}(a)| = 2, \\ 
m=4: \quad 
\chi(\calC,t)&=(t-1)(t-11)(t-13)(t-15), \quad 
|\Ch(\calC)| = 5376, \quad 
|\varphi_{\calA}^{-1}(a)| = 14, \\ 
m=5: \quad 
\chi(\calC,t)&=(t-1)(t-29)(t-31)(t^2-60t+971), \\ 
& |\Ch(\calC)| = 1981440, \quad 
|\varphi_{\calA}^{-1}(a)| = 516, \\ 
m=6: \quad 
\chi(\calC,t)&=(t-1)(t^5-363t^4+54310t^3-4182690t^2+165591769t-2691439347), \\ 
& |\Ch(\calC)| = 5722536960, \quad 
|\varphi_{\calA}^{-1}(a)| = 124187. 
\end{align*}

Note that the characteristic polynomial $\chi(\calC, t)$ factors 
into polynomials of degree one over $\bbZ$ for $m \le 4$.

%

\ 

\begin{flushleft}
{\bf Acknowledgments} 
\end{flushleft} 
The authors are very grateful to an anonymous referee 
for 
valuable comments on an earlier version of this paper. 
The example of the Catalan arrangement in 
Subsection \ref{subsec:Catalan} was suggested by 
the referee.

\end{document}